\documentclass[12pt]{article}

\title{\textbf{Two philosophical applications of algorithmic
information theory}}

\author{Gregory Chaitin\thanks{IBM Research Division, P. O. Box 218,
Yorktown Heights, NY 10598, USA. E-mail: chaitin@us.ibm.com.}
}

\date{}

\begin{document}
\sloppy

\maketitle

\begin{abstract}
Two philosophical applications of the concept of program-size complexity
are discussed. First, we consider the light program-size complexity
sheds on whether mathematics 
is invented or discovered, i.e., is empirical or is \emph{a priori.}
Second, we propose that the notion of algorithmic independence 
sheds light on the question of being and how the world of our experience
can be partitioned into separate entities.
\end{abstract}

\section*{1. Introduction. Why is program size of philosophical interest?}

The cover of the January 2003 issue of \emph{La Recherche} asks this dramatic question:
\begin{center}
Dieu est-il un ordinateur? [Is God a computer?]
\end{center}
The long cover story [1] is a reaction to Stephen Wolfram's controversial
book \emph{A New Kind of Science} [2].  The first half of the article 
points out Wolfram's predecessors, and the second half criticizes Wolfram.

The second half of the article begins (p.\ 38) with these words:

\begin{quote}
Il [Wolfram] n'avance aucune raison s\'erieuse de
penser que les complexit\'es de la nature
puissent \^etre g\'en\'er\'ees par des r\`egles
\'enon\c{c}ables sous forme de programmes
informatiques simples.
\end{quote}

The reason for thinking that a simple program
might describe the world is, basically, just Plato's
postulate that the universe is rationally
comprehensible (\emph{Timaeus}).  A sharper statement
of this principle is in Leibniz's \emph{Discours de
m\'etaphysique} [3], section \textbf{VI}. Here is Leibniz's original French (1686):

\begin{quote}
Mais Dieu a choisi celuy
qui est le plus parfait, c'est \`a dire celuy qui est
en m\^{e}me temps le plus simple en hypotheses
et le plus riche en phenomenes, comme pourroit
estre une ligne de Geometrie dont la construction
seroit ais\'ee et les propriet\'es et effects seroient
fort admirables et d'une grande \'etendue.
\end{quote}

For an English translation of this, see [4].

And Hermann Weyl [5] discovered that 
in \emph{Discours de m\'etaphysique} 
Leibniz also states that
a physical law has no explicative power if it is
as complicated as the body of data it was invented to explain.\footnote
{See the Leibniz quote in Section 2.3 below.}  

This is where algorithmic information theory (AIT) comes in.
AIT posits that a theory that explains $X$ is a computer program for calculating $X$, 
that therefore
must be smaller, much smaller, than the size in bits of the data $X$ that it explains.
AIT makes a decisive contribution to philosophy 
by providing a mathematical theory of complexity.
AIT defines the \emph{complexity} or \emph{algorithmic information content} of $X$
to be the size in bits $H(X)$ of the smallest computer program for calculating $X$.
$H(X)$ is also the complexity of the most elegant (the simplest) theory for $X$.

In this article we discuss some other philosophical applications of AIT.

\bigskip

For those with absolutely no background in philosophy,
let me recommend two excellent introductions,
Magee [6] and Brown [7].  
For introductions to AIT, see Chaitin [8, 9].
For another discussion of the philosophical implications of AIT,
see Chaitin [10].

\section*{2. Is mathematics empirical or is it a priori?}

\subsection*{2.1. Einstein: math is empirical}

Einstein was a physicist and he believed that math is invented, not discovered.
His sharpest statement on this is his declaration that
``the series of integers is obviously
an invention of the human mind, a self-created tool
which simplifies the ordering of certain sensory experiences.''

Here is more of the context:

\begin{quote}
In the evolution of philosophic thought through the centuries
the following question has played a major r\^ole: What knowledge
is pure thought able to supply independently of sense perception?
Is there any such knowledge?\ldots\ I am convinced that\ldots\
the concepts which arise in our thought and in our linguistic
expressions are all\ldots\ the free creations of thought which
can not inductively be gained from sense-experiences\ldots\
\textbf{Thus, for example, the series of integers is obviously an
invention of the human mind, a self-created tool which simplifies
the ordering of certain sensory experiences}.\footnote
{[The boldface emphasis in this and future quotations is mine,
not the author's.]}
\end{quote}

The source is Einstein's essay
``Remarks on Bertrand Russell's theory
of knowledge.''  It was published in 1944 in the volume [11] on
\emph{The Philosophy of Bertrand Russell} edited by
Paul Arthur Schilpp, and it was reprinted in 1954 in Einstein's
\emph{Ideas and Opinions} [12].

And in his \emph{Autobiographical Notes} [13] Einstein
repeats the main point of his Bertrand Russell essay, in
a paragraph on Hume and Kant in which he
states that ``all concepts, even those closest to
experience, are from the point of view of logic freely chosen
posits.''
Here is the bulk of this paragraph:

\begin{quote}
Hume saw clearly that certain concepts, as for example that of causality,
cannot be deduced from the material of experience by logical methods.
Kant, thoroughly convinced of the indispensability of certain concepts,
took them\ldots\ to be the necessary premises of any kind of thinking
and distinguished them from concepts of empirical origin. I am convinced,
however, that this distinction is erroneous or, at any rate, that it
does not do justice to the problem in a natural way. \textbf{All concepts,
even those closest to experience, are from the point of view of logic
freely chosen posits}\ldots
\end{quote}

\subsection*{2.2. G\"odel: math is a priori}

On the other hand, G\"odel was a Platonist and believed that math is a priori.
He makes his position blindingly clear in the introduction to an unpublished lecture
G\"odel *1961/?, ``The modern development of the foundations
of mathematics in the light of philosophy,''
\emph{Collected Works} [14], vol.\ 3:\footnote
{The numbering scheme used in G\"odel's \emph{Collected Works} 
begins with an * for unpublished papers,
followed by the year of publication, or the first/last year that G\"odel
worked on an unpublished paper.}

\begin{quote}
I would like to attempt here to describe, in terms of philosophical
concepts, the development of foundational research in mathematics\ldots,
and to fit it into a general schema
of possible philosophical world-views [Weltanschauungen]\ldots\
I believe that the most fruitful principle for gaining an overall
view of the possible world-views will be to divide them up according
to the degree and the manner of their affinity to or, respectively,
turning away from metaphysics (or religion). In this way we
immediately obtain a division into two groups: skepticism, materialism
and positivism stand on one side, spiritualism, idealism and theology
on the other\ldots\ Thus one would, for example, say that apriorism
belongs in principle on the right and empiricism on the left side\ldots\
Now it is a familiar fact, even a platitude, that the development
of philosophy since the Renaissance has by and large gone from right
to left\ldots\ It would truly be a miracle if this (I would like to
say rabid) development had not also begun to make itself felt in the
conception of mathematics. Actually, \textbf{mathematics, by its nature as
an a priori science}, always has, in and of itself, an inclination
toward the right, and, for this reason, \textbf{has long withstood the spirit
of the time} [Zeitgeist] that has ruled since the Renaissance; i.e.,
the empiricist theory of mathematics, such as the one set forth by
Mill, did not find much support\ldots\
Finally, however, around the turn of the century, its hour struck: in
particular, it was the antinomies of set theory, contradictions that 
allegedly appeared within mathematics, whose significance was 
exaggerated by skeptics and empiricists and which were employed as
a pretext for the leftward upheaval\ldots\
\end{quote}

Nevertheless, the Platonist G\"odel makes some remarkably strong statements
in favor of adding to mathematics axioms which are not self-evident and which
are only justified pragmatically.
What arguments does he present in support of these heretical views?
 
First let's take a look at his discussion of whether
Cantor's continuum hypothesis could be established using a new axiom [G\"odel 1947, 
``What is Cantor's continuum problem?'', \emph{Collected Works,} vol.\ 2]:

\begin{quote}
\ldots even \textbf{disregarding the intrinsic necessity of some new axiom}, 
and even in case it has no intrinsic
necessity at all, \textbf{a probable decision about its truth is possible also in another way, namely,
inductively} by studying its ``success.'' Success here means fruitfulness in consequences, in
particular in ``verifiable'' consequences, i.e., consequences demonstrable without the new axiom,
whose proofs with the help of the new axiom, however, are considerably simpler and easier to
discover, and make it possible to contract into one proof many different proofs. The axioms for
the system of real numbers, rejected by intuitionists, have in this sense been verified to some
extent, owing to the fact that analytical number theory frequently allows one to prove
number-theoretical theorems which, in a more cumbersome way, can subsequently be verified by
elementary methods. A much higher degree of verification than that, however, is conceivable.
\textbf{There might exist axioms} so abundant in their verifiable consequences, shedding so much light
upon a whole field, and yielding such powerful methods for solving problems (and even solving
them constructively, as far as that is possible) \textbf{that}, no matter whether or not they are intrinsically
necessary, they \textbf{would have to be accepted at least in the same sense as any well-established
physical theory}. 
\end{quote}

Later in the same paper G\"odel restates this: 

\begin{quote}
It was pointed out earlier\ldots\ that, \textbf{besides mathematical intuition, 
there exists another} (though only probable) 
\textbf{criterion of the truth of mathematical axioms, namely} their \textbf{fruitfulness} in mathematics
and, one may add, possibly also in physics\ldots\ 
The simplest case of an application of the criterion
under discussion arises when some\ldots\ axiom has number-theoretical consequences verifiable by
computation up to any given integer. 
\end{quote}

And here is an excerpt from G\"odel's contribution  
[G\"odel 1944, ``Russell's mathematical logic,'' \emph{Collected Works,} vol.\ 2]
to the same Bertrand Russell festschrift volume [11] that was quoted above:

\begin{quote}
The analogy between mathematics and a natural science is enlarged upon by Russell also in
another respect\ldots\ \textbf{axioms need not be evident in themselves, but rather their justification lies
(exactly as in physics) in the fact that they make it possible for these ``sense perceptions'' to be
deduced}\ldots\ 
I think that\ldots\ this view has been largely justified by subsequent developments, and it is
to be expected that it will be still more so in the future. It has turned out that the solution of
certain arithmetical problems requires the use of assumptions essentially transcending
arithmetic\ldots\ 
Furthermore it seems likely that for deciding certain questions of abstract set theory
and even for certain related questions of the theory of real numbers new axioms based on some
hitherto unknown idea will be necessary. Perhaps also the apparently insurmountable difficulties
which some other mathematical problems have been presenting for many years are due to the fact
that the necessary axioms have not yet been found. Of course, under these circumstances
mathematics may lose a good deal of its ``absolute certainty;'' but, under the influence of the
modern criticism of the foundations, this has already happened to a large extent\ldots 
\end{quote}

Finally, take a look at this excerpt 
from G\"odel *1951, ``Some basic theorems on the foundations,''
\emph{Collected Works,} vol.\ 3,
an unpublished essay by G\"odel:

\begin{quote}
I wish to point out that one may conjecture the truth of a universal proposition
(for example, that I shall be able to verify a certain property for \emph{any}
integer given to me) and at the same time conjecture that no general proof for this fact
exists. It is easy to imagine situations in which both these conjectures would be very
well founded. For the first half of it, this would, for example, be the case if the
proposition in question were some equation $F(n)=G(n)$ of two number-theoretical
functions which could be verified up to \emph{very} great numbers $n$.\footnote
{Such a verification of an \emph{equality} (not an inequality) between two number-theoretical
\textbf{functions of not too complicated or artificial structure} would certainly give a
great probability to their complete equality, although its numerical value could
not be estimated in the present state of science. However, it is easy to give
examples of general propositions about integers where the probability can be estimated 
even now\ldots}
Moreover, exactly as in the natural sciences, this \emph{inductio per enumerationem simplicem}
is by no means the only inductive method conceivable in mathematics. 
I admit that every mathematician has an inborn abhorrence to giving more than
heuristic significance to such inductive arguments.
I think, however,
that this is due to the very prejudice that mathematical objects somehow have no
real existence. \textbf{If mathematics describes an objective
world just like physics, there is no reason why inductive methods should not be
applied in mathematics just the same as in physics}. The fact is that in mathematics
we still have the same attitude today that in former times one had toward all science,
namely, we try to derive everything by cogent proofs from the definitions (that is,
in ontological terminology, from the essences of things). Perhaps this method, if it
claims monopoly, is as wrong in mathematics as it was in physics.
\end{quote}

So G\"odel the Platonist has nevertheless managed to arrive, at least partially, at what I would  
characterize, following Tymoczko [16], as a pseudo-empirical or a quasi-empirical position!

\subsection*{2.3. AIT: math is quasi-empirical}

What does algorithmic information theory have to contribute to this discussion?
Well, I believe that AIT also supports a quasi-empirical view of mathematics.
And I believe that it provides further justification for G\"odel's belief
that we should be willing to add new axioms.

Why do I say this?

As I have argued on many occasions, AIT, by measuring the complexity 
(algorithmic information content) of axioms
and showing that G\"odel incompleteness is natural and ubiquitous, deepens the arguments
that forced G\"odel, in spite of himself, in spite of his deepest instincts about the nature
of mathematics, to believe in inductive mathematics.
And if one considers the use of induction rather than deduction to establish 
mathematical facts, some kind of notion of \emph{complexity} must necessarily be
involved. For as Leibniz stated in 1686, a theory is
only convincing to the extent that it is substantially \emph{simpler} than the
facts it attempts to explain:

\begin{quote}
\ldots non seulement rien n'arrive dans le monde, qui soit absolument irregulier, mais
on ne s\c{c}auroit m\^emes rien feindre de tel. Car supposons par exemple que quelcun
fasse quantit\'e de points sur le papier \`a tout hazard, comme font ceux qui exercent
l'art ridicule de la Geomance, je dis qu'il est possible de trouver une ligne
geometrique dont la motion soit constante et uniforme suivant une certaine regle, en
sorte que cette ligne passe par tous ces points\ldots\ Mais \textbf{quand une regle est
fort compos\'ee, ce qui luy est conforme, passe pour irr\'egulier}. Ainsi on peut
dire que de quelque maniere que Dieu auroit cr\'e\'e le monde, il auroit
tousjours est\'e regulier et dans un certain ordre general. Mais Dieu a choisi
celuy qui est le plus parfait, c'est \`a dire celuy qui est en m\^eme temps \textbf{le plus
simple en hypotheses et le plus riche en phenomenes}\ldots\ 
[\emph{Discours de m\'etaphysique,} \textbf{VI}]
\end{quote}

In fact G\"odel himself, in considering inductive rather than deductive mathematical proofs, 
began to make some tentative initial attempts to formulate and utilize notions of 
complexity. (I'll tell you more about this in a moment.)  And it is here that AIT
makes its decisive contribution to philosophy, by providing a highly-developed and
elegant mathematical theory of complexity. How does AIT do this? It does this by considering 
the size of the smallest computer
program required to calculate a given object $X$, which may also be considered to be
the most elegant theory that explains $X$.

Where does G\"odel begin to think about complexity?
He does so in two footnotes in vol.\ 3 of his \emph{Collected Works.}
The first of these is a footnote to G\"odel *1951. 
This footnote begins ``Such a verification\ldots''\@ and it was reproduced, 
in part, in Section 2.2 above.
And here is the relevant portion of the second, the more interesting, 
of these two footnotes: 

\begin{quote}
\ldots Moreover, if every number-theoretical question of Goldbach type\ldots\
is decidable by a mathematical proof, there \emph{must} exist an infinite set of
independent evident axioms, i.e., a set $m$ of evident axioms which are not
derivable from \emph{any} finite set of axioms (no matter whether or not the latter
axioms belong to $m$ and whether or not they are evident). Even if solutions are
desired only for all those problems of Goldbach type which are simple enough to be 
formulated in a few pages, \textbf{there must exist} a great number of evident axioms or
evident \textbf{axioms of great complication, in contradistinction to the few simple axioms
upon which all of present day mathematics is built}. (It can be proved that, in order
to solve all problems of Goldbach type of a certain degree of complication $k$, one
needs a system of axioms whose degree of complication, up to a minor correction, is
$\ge k$.)\footnote
{[This is reminiscent of the theorem in AIT that $p_k =$ (the program of size $\le k$
bits that takes longest to halt) is the simplest possible ``axiom'' from which one can solve the 
halting problem for all programs of size $\le k$.
Furthermore, $p_k$'s size and complexity
both differ from $k$ by at most a fixed number of bits: $|p_k|=k+O(1)$ and
$H(p_k)=k+O(1)$. 

Actually, 
in order to solve the halting problem for all programs of size $\le k$,
in addition to $p_k$ one needs to know $k-|p_k|$,
which is how much $p_k$'s size differs from $k$. 
This fixed amount of additional information is required in order to be able to
determine $k$ from $p_k$.]}
\end{quote}

This is taken from G\"odel *1953/9--III, one of the versions of his
unfinished paper ``Is mathematics syntax of language?''\@
that was intended for, but was finally not included, in Schilpp's Carnap festschrift 
in the same series as the Bertrand Russell festschrift [11].

Unfortunately these tantalizing glimpses are, as far as I'm aware, all that we
know about G\"odel's thoughts on complexity. Perhaps volumes 4 and 5, the two final volumes of
G\"odel's \emph{Collected Works,} which contain G\"odel's correspondence with other mathematicians, 
and which will soon be available, will shed further light on this.

\bigskip 

Now let me turn to a completely different---but I believe equally fundamental---application of AIT.

\section*{3. How can we partition the world into distinct entities?}

For many years I have asked myself, ``What is a living being? How can we define
this mathematically?!'' I still don't know the answer!  But at least I think I now
know how to come to grips with the more general notion of ``entity'' or ``being.''
In other words, how can we decompose our experience into parts?  How can we
partition the world into its components?   By what right do we do this in spite of
mystics who like Parmenides insist that the world must be perceived as an organic unity
(is a single substance)
and \textbf{cannot} be decomposed or analized into independent parts?

I believe that the key to answering this fundamental question lies in AIT's
concept of \emph{algorithmic independence.}  What is algorithmic independence?
Two objects $X$ and $Y$ are said to be algorithmically independent if their
complexity is (approximately) additive.  In other words, $X$ and $Y$ are
algorithmically independent if their information content decomposes additively,
i.e., if their joint information content (the information content of $X$ \textbf{and} $Y$)
is approximately equal to the sum of their individual information contents:
\[
    H(X,Y) \approx H(X) + H(Y) .
\]
More precisely, the left-hand side is the size in bits of the smallest program
that calculates the pair $X$, $Y$, and the right-hand side adds the size in bits
of the smallest program that produces $X$ to the size in bits of 
the smallest program that calculates $Y$.

Contrariwise, if $X$ and $Y$ are \textbf{not at all} independent, then 
it is much better to compute them together than to compute them separately
and $H(X)+H(Y)$ will be much larger than $H(X,Y)$.
The worst case is $X=Y$. Then $H(X)+H(Y)$ is twice as large as $H(X,Y)$.

I feel that this notion of algorithmic independence is the key to decomposing the world into parts,
parts the most interesting example of which are living beings, particularly human beings.
For what enables me to partition the world in this way?  The fact that thinking of the
world as a sum of such parts does not complicate my description of the world substantially
and at the same time enables me to use separate subroutines such as ``my wife'' and ``my cat''
in thinking about the world.  That is why such an analysis of the world, such a
decomposition, works.   

Whereas on the contrary ``my left foot'' and ``my right hand''
are not well thought of as independent components of the world but can best be understood 
as parts of me.   A description of my right hand and its activities and history would
not be substantially simpler than a description of me and my entire life history,
since my right hand is a part of me whose actions express my intentions, and not its own
independent desires.
 
Of course, these observations are just the beginning. A great deal more work
is needed to develop this point of view\ldots

\bigskip

For a technical discussion of algorithmic independence and
the associated notion of
\emph{mutual algorithmic information} defined as follows
\[
H(X:Y) \equiv H(X)+H(Y)-H(X,Y),
\]
see my book Chaitin [17].

\section*{4. Conclusion and future prospects}

Let's return to our starting point, to the cover of the January 2003 issue of
\emph{La Recherche}.  Is God a computer, as Wolfram and some others think,
or is God, as Plato and Pythagoras affirm, a mathematician?

And, an important part of this question,
\textbf{is the physical universe discrete}, the way computers prefer, \textbf{not continuous}, 
the way it seems to be in classical Newtonian/Maxwellian physics?
Speaking personally,
I like the discrete, not the continuous. 
And my theory, AIT, deals with discrete, digital information, bits, not with continuous 
quantities.
But the physical universe
is of course free to do as it likes!

Hopefully pure thought will not be called upon to resolve this.
Indeed, I believe that it is incapable of doing so; Nature will have
to tell us. Perhaps someday an \emph{experimentum crucis} will 
provide a definitive answer.  
In fact, for a hundred years quantum physics has been pointing
insistently in the direction of discreteness.\footnote
{Discreteness in physics actually began even earlier,  
with atoms.  And then, my colleague John Smolin points out,
when Boltzmann introduced coarse-graining in statistical mechanics.}

\section*{References}

\begin{itemize}
\item[{[1]}]
O. Postel-Vinay,
``L'Univers est-il un calculateur?''\@ [Is the universe a calculator?],
\emph{La Recherche,} no.\ 360, January 2003, pp.\ 33--44.
\item[{[2]}]
S. Wolfram, \emph{A New Kind of Science,} Wolfram Media, 2002.
\item[{[3]}]
Leibniz, \emph{Discours de m\'etaphysique,} Gallimard, 1995.
\item[{[4]}]
G. W. Leibniz, \emph{Philosophical Essays,} Hackett, 1989.
\item[{[5]}]
H. Weyl, \emph{The Open World,} Yale University Press, 1932, Ox Bow Press, 1989.
\item[{[6]}]
B. Magee, \emph{Confessions of a Philosopher,} Modern Library, 1999.
\item[{[7]}]
J. R. Brown, \emph{Philosophy of Mathematics,} Routledge, 1999.
\item[{[8]}]
G. J. Chaitin, ``Paradoxes of randomness,'' \emph{Complexity,} vol.\ 7, no.\ 5, pp.\ 14--21, 2002.
\item[{[9]}]
G. J. Chaitin, ``Meta-mathematics and the foundations of mathematics,''
\emph{Bulletin EATCS,} vol.\ 77, pp.\ 167--179, 2002.
\item[{[10]}]
G. J. Chaitin, ``On the intelligibility of the universe and the notions of
simplicity, complexity and irreducibility,''
\\
http://arxiv.org/abs/math.HO/0210035, 2002.
\item[{[11]}]
P. A. Schilpp, \emph{The Philosophy of Bertrand Russell,} Open Court, 1944.
\item[{[12]}]
A. Einstein, \emph{Ideas and Opinions,} Crown, 1954, Modern Library, 1994.
\item[{[13]}]
A. Einstein, \emph{Autobiographical Notes,} Open Court, 1979.
\item[{[14]}]
K. G\"odel, \emph{Collected Works,} vols.\ 1--5, Oxford University Press, 1986--2003.
\item[{[15]}]
\emph{Kurt G\"odel: Wahrheit \& Beweisbarkeit} [truth and provability], vols.\ 1--2, \"obv \& hpt, 2002.
\item[{[16]}]
T. Tymoczko, \emph{New Directions in the Philosophy of Mathematics,}
Princeton University Press, 1998.
\item[{[17]}]
G. J. Chaitin, \emph{Exploring Randomness,} Springer-Verlag, 2001.
\end{itemize}
Chaitin's papers are also available at 
\\
http://cs.auckland.ac.nz/CDMTCS/chaitin.

\end{document}